\documentstyle[12pt]{article}
 \newcommand{\beq}{\begin{equation}}
\newcommand{\eeq}{\end{equation}}
 \newcommand{\beth}{\begin{theo}}
\newcommand{\eth}{\end{theo}}

\newcommand{\Subset}{\subset\subset}

\newcommand{\Ifo}{{ \Psi_{f,0}}}
\newcommand{\Ifxm}{{\Psi_{f,x^m}}}
\newcommand{\Ifxo}{{\Psi_{f,x^0}}}
\newcommand{\Iuxm}{{\Psi_{v,x^m}}}
\newcommand{\Ivxm}{{\Psi_{v,x^m}}}

\newcommand{\tfxo}{\tau_{f}(x^0)}
\newcommand{\tfo}{\tau_{f}(0)}
\newcommand{\obl}{\Omega}
\newcommand{\nol}{\{0\}}
\newcommand{\okn}{1\le k\le n}
\newcommand{\ojn}{1\le m\le N}
\newcommand{\Ij}{\Psi_m}
\newcommand{\hpsi}{\Psi}
\newcommand{\green}{G_{\Phi,\obl}}
\newcommand{\mera}{T_{\Phi }}
\newcommand{\gree}{G}
\newcommand{\ifoa}{I(F,0,a)}
\newcommand{\ifxa}{I(F,x^0,a)}
\newcommand{\nfoa}{n(f,0,a)}
\newcommand{\nfcoa}{n(f_c,0,a)}
\newcommand{\nfxa}{n(f,x^0,a)}

\newcommand{\ntxo}{\nu(T,x^0)}
\newcommand{\vaa}{\varphi_a}
\newcommand{\cldo}{ {D_1}}
\newcommand{\cld}{ \overline{D}}
\newcommand{\dvapi}{}

\newcommand{\Rn}{{ \bf R}^n}
\newcommand{\Rnm}{{ \bf R}_-^n}
\newcommand{\Rm}{{ \bf R}_-}
\newcommand{\Cn}{{ \bf C}^n}
\newcommand{\Cp}{{ \bf C}^p}
\newcommand{\Nn}{{ \bf N}^n}

\newtheorem{theo}{Theorem}
\newtheorem{cor}{Corollary}

\newtheorem{prop}{Proposition}

\begin{document}

\begin{center}
 {\bf \large  LOCAL INDICATORS FOR}
\end{center}

\begin{center}
{\bf \large  PLURISUBHARMONIC FUNCTIONS}
\end{center}

\vskip0.5cm
\begin{center}
{\bf By Pierre LELONG and Alexander RASHKOVSKII  }
\end{center}

\vskip1cm

  {\small {\sc Abstract.} - The notion of index, classical in number
theory and extended in \cite{keyL4} to plurisubharmonic functions,
allows to define an
indicator which  is applied to the study of the Monge-Amp\`ere
operator and   a pluricomplex Green function.  }


\vskip 0.5cm

\section{Introduction}
We recall some local notions which are often used in various
investigations.

\medskip
$a)$.
The  index $\ifxa$
of a zero $x^0\in\obl$ of a holomorphic function $F\in Hol(\obl)$, where
$\obl$ is a domain of $\Cn$,  is used in important results of number theory
\cite{keyLa}; $\ifxa$ is defined by
means of the set $\omega\in\Nn$ of $n$-tuples $(i)=(i_1,\ldots,i_n)$ such
that
$ D^{(i)}F(x^0)\neq 0$. Given a direction
$(a)=(a_k>0; \: \okn)\in\Rn$,
 \beq
\ifxa=\inf_{(i)}(a,i) \quad {\rm for}\ (a,i)= \sum_ka_ki_k\ge 0 \ {\rm and}\
(i)\in\omega.
\label{eq:L1}
\eeq

\medskip
$b)$. In fact, the index $\ifxa$ is a property (function of $a$ and $x^0$)
of the current of integration $[W]= \dvapi dd^c\log|F|$ over the analytic
set $W=\{x\in\obl:\: F(x)=0\}$,
where $d=\partial + \bar\partial$ and $d^c= (2\pi i)^{-1}( \partial
-\bar\partial)$.
 We denote by $PSH(\obl)$ the class of
plurisubharmonic
 functions in a domain
$\obl$ of $\Cn$ and by $\Theta_p(\obl)$ the class of positive closed
currents represented by homogeneous forms of $dx_k,\:d\bar x_k,
\ \okn$, of bidegree $(n-p,n-p)$. The Lelong number $\ntxo$ at $x^0\in\obl$
for $T\in \Theta_p(\obl)$ is related to the trace measure of $T$,
\beq
\sigma_T=T\wedge\beta_p,
\label{eq:L2}
\eeq
where $\beta_p=(p!)^{-1}\beta_1^p$ is the volume element of $\Cp$.
By the definition,
\beq
\ntxo=\lim_{r\to 0}\,(\tau_{2p}r^{2p})^{-1}\sigma_T[B^{2n}(x^0,r)]
\label{eq:L3}
\eeq
where $\tau_{2p}$ is the volume of the unit ball $B^{2p}(0,1)$ of $\Cp$.
In (\ref{eq:L3}), the trace measure $\sigma_T$ belongs to the remarkable
class ${\cal P}_p(\obl)$ of positive measures characterized by the property
that the quotient in the right hand side of (\ref{eq:L3}) is an increasing
function of $r$ (see \cite{keyL5}). Another definition of the number $\ntxo$
is derived from (\ref{eq:L2}) and (\ref{eq:L3}) by setting
$\varphi_1(x)=\log|x-x^0|$:
\beq
\ntxo=\lim_{t\to -\infty}\,\int_{\{\varphi_1<t\}}T\wedge  (dd^c
\varphi_1)^{p}.
\label{eq:L4}
\eeq

\medskip
$c)$. By replacing in (\ref{eq:L4}) $\varphi_1$ with a function $\varphi\in
PSH(\obl)$ such that $\exp\varphi$ is continuous and the set
$\{\varphi(x)=-\infty\}$ is relatively compact, J.-P.~Demailly (see
\cite{keyD1} and \cite{keyD2}) has found new  important applications of
(\ref{eq:L4}). The choice of the function
\beq
\varphi_a(x)=\sup_k a_k^{-1}\log|x_k-x_k^0|
\label{eq:L5}
\eeq
which is circled with the center $x^0$ and such that
$\{\varphi_a(x)=-\infty\}$ is reduced to $x^0$, allows us to put
into this framework the notion of index $\ifxa$ of a zero of a
holomorphic function $F$
at $x^0$. As was proved in \cite{keyL4},
\beq
\ifxa=(a_1\cdot a_2\ldots a_n)\,\nu[T,x^0,\varphi_a(x-x^0)]
\label{eq:L6}
\eeq
for $T= \dvapi dd^c\log|F|\in\Theta_1(\obl)$. It can be extended to
arbitrary plurisubharmonic functions by setting for
  $f\in PSH(\obl)$ and
$x^0\in\obl$,
\beq
\nfxa=(a_1\cdot a_2\ldots a_n)\,\nu[ \dvapi dd^cf,x^0,\varphi_a(x-x^0)],
\label{eq:L7}
\eeq
so that
\beq
\ifxa=\nfxa \quad {\rm for}\ f=\log|F|.
\label{eq:L8}
\eeq

\medskip
$d)$. In the recent paper \cite{keyL4}, $\nfxa$ has appeared in a simpler
form given for $x^0=0$ by the relation
\beq
\nfoa=\lim_{w\to 0}\: (\log{ w})^{-1} f(w,x,a),\quad 0\le w\le 1,\
x\in\obl\setminus A,
\label{eq:L9}
\eeq
where $f(w,x,a)=f(w^{a_1}x_1,\ldots,w^{a_n}x_n)$ and $A$ is an algebric set
for $f=\log|F|$. In the general case, for $f\in PSH(\obl)$, $A$ is of
zero measure and
$$
\nfoa=\liminf_{w\to 0}\:
(\log{ w})^{-1}f(w,x,a)
$$
 outside  a set $A'\subset A$
which is pluripolar in $\obl$.

\medskip
In the first part of the present paper we will use (\ref{eq:L9}) for a
study of singularities of plurisubharmonic functions $f$,
 supposing
$f\in PSH(\obl)$ and
$D\Subset\obl$, where $D$ is the   unit polydisk
$\{x\in\Cn:\:\sup |x_k|< 1\}$; we denote by $PSH_-(D)$
the class of $f$ satisfying $sup\, \{f(x): x\in D\}\le 0$   and
$f\not\equiv -\infty$.
The domain $D$ as well as the weights $\vaa$
in (\ref{eq:L5})
being circled, it is natural to work with a circled image $f_c$ of $f$ and
then with  its convex image   on the space $\Rnm$ of $u_k=\log|x_k|,\
\okn$.
Developing the results of J.-P.~Demailly and C.O.~Kiselman \cite{keyK2}
 we get the value $\nfoa$ which produces, via $a_k=-\log|y_k|$,
a function $\Ifo(y) $, the local indicator of $f$ at $0$,
which is plurisubharmonic in
$D_{(y)}$, the unit polydisk in the space $\Cn_{y}$. It
satisfies the
 Monge-Amp\`ere equation
\beq
(dd^c\Ifo)^n=\tfo\,\delta(0),\quad \tfo>0,
\label{eq:L10}
\eeq
where $\delta(0)$ is the Dirac measure at the origin of $\Cn_y$, and
\beq
(dd^cf)^n\ge\tfxo\,\delta(x^0)
\label{eq:L11}
\eeq
 for $f\in PSH(\obl)$ such that $(dd^cf)^n$ is well defined, and
$x^0\in\obl$. Moreover, there is the relation
$\tfxo\ge [\nu(T,x^0)]^n$ for $T=dd^cf\in\Theta_{n-1}
(\obl).$

Then we consider a class of plurisubharmonic functions on $\obl$ with
singularities on
a finite set $\{x^1,\ldots,x^N\}$, controlled by given indicators
$\Ij,\ \ojn$. We construct a plurisubharmonic function $\gree$
vanishing on $\partial\obl$ and such that
$\Psi_{ \gree,x^m}=\Ij,\,\ojn$,
and $(dd^c G)^n=\sum_m\,\tau_m\delta(x^m)$ with $\tau_m$ the mass
of $(dd^c\Ij)^n$. We prove that $G$ is the unique plurisubharmonic function
with these properties. For the case $\Ij(x)=\nu_m\log|x-x^m|$, it
coincides with the pluricomplex Green function with weighted poles at
$x^1,\ldots,x^N$ (\cite{keyL3}).
 We prove a variant of comparison theorem for plurisubharmonic functions
with controlled singularities and study a Dirichlet problem for this class
of functions.

A part of the results of Section 3 are close to those from \cite{keyZ}
where a weighted pluricomplex Green function with infinite singular set was
introduced and the corresponding Dirichlet problem was studied.

\section{ Circled functions   and convex projections}

We will consider here the case $x^0=0$ and $ f\in PSH(\obl)$   supposing
$0\in D\Subset\obl$, where $D$ is the   unit polydisk
$\{x\in\Cn:\:\sup |x_k|< 1\}$, and $f\in PSH_-(D)$, that is
 $f(x)\le 0$ for $x\in\cld$ and $f\not\equiv
-\infty$.

A set $A\subset\Cn$ is  called {\it $0$-circled} (or just {\it circled})
if $x=(x_k)\in A$ implies $x'=(x_ke^{i\theta_k})\in A$ for $0\le\theta_k
\le 2\pi,\ \okn$.
We will say that a function $f(x)$ defined on $A$, is
circled if
it is invariant with respect to the rotations $x_k\mapsto x_ke^{i\theta_k},
\ \okn$.

Given a function $f\in PSH(\obl),\ \obl$ being a circled domain,  we
 consider a circled function  $f_c\in PSH(\obl)$ equal to the mean value of
$f(x_ke^{i\theta_k})$ with respect to $0\le\theta_k
\le 2\pi,\ \okn$. In what follows, we will also use another circled
function $f_c'\ge f_c$, equal to the maximum of $f(x_ke^{i\theta_k})$
for $0\le\theta_k\le 2\pi,\ \okn$.
Note that the differential operators,
namely $\partial, \ \bar\partial,\ d=\partial + \bar\partial$ and
$d^c= (2\pi i)^{-1}( \partial
-\bar\partial)$, commute with the mapping $f\mapsto f_c$,
so $(\partial f)_c=\partial f_c$,  however it is not the case for
$f\mapsto f'_c$.

To a Radon measure $\sigma$ on a circled domain $\obl$, we  relate
a circled measure $\sigma_c$ defined by $\sigma_c(f)=\sigma(f_c)$
for continuous functions $f$. In the
same way, to a current $T\in
\Theta_p(\obl)$ we  associate a circled current $T_c$ which is defined
on homogeneous forms $ \lambda$ of bidegree
$(p,p)$ by $T_c( \lambda)=T( \lambda_c)$, where $\lambda_c$ are obtained by
replacing the coefficients of $ \lambda$ with their mean values
with respect to $\theta_k$.
 In particular, if $T=dd^cf,\ f\in PSH(\obl)$
and $\obl$ circled, $T_c=dd^cf_c$. It gives us a specific property
of the value $\nfoa$ defined by (\ref{eq:L4}) and (\ref{eq:L6}), and of
the index $\ifoa$.

\begin{prop}
Let $\obl$ be a $0$-circled domain and $T\in
\Theta_p(\obl)$.  For every $0$-circled weight $\varphi$,
$\nu(T,\varphi)=\nu(T_c,\varphi)$.
In particular, since the weight $\vaa$ defined by (\ref{eq:L5}) for $x^0=0$,
 is circled,
  the number $\nfoa$  in
(\ref{eq:L7}) for $f\in PSH(\obl)$ at the origin can be calculated by
replacing $f$ with $f_c$:
 $\nfoa=\nfcoa$, and the number $\nfoa$ is calculated on the convex
image $g(u)$ of $f_c$, $g(u)=f[\exp (u_k+i\theta_k)]$:
$$
\nfoa=\lim_{v\to -\infty}\,v^{-1}g(u_k+a_kv).
$$
\label{prop:L1}
\end{prop}

 \medskip
 For $f\in PSH_-(D)$, as was shown in  \cite{keyL4},
 $$
\nfoa=\lim_{w\to 0}\,(\log w)^{-1}f(w^{a_1}x_1,\ldots,w^{a_n}x_n)
$$
for almost all $x\in D$. The limit exists for all $x\in D$ when replacing
$f(x)$ by $f_c(x)$ or by $f_c'(x)$. Indeed, for $R>1$
\beq
f_c(x)\le f_c'(x)\le\gamma_R\,f_c(Rx)\le 0
\label{eq:L11}
\eeq
with $\gamma_R=(R-1)^n(R+1)^{-n}$ which satisfies $1-\epsilon\le\gamma_R\le
1$ for $R>R_0(\epsilon)$.

\medskip

 The calculation of $\nfoa$ for
$f\in PSH_-(D)$  uses the convex image
$g_f(u)=f_c[\exp (u_k+i\theta_k)]$ or
$g_f'(u)=f_c'[\exp (u_k+i\theta_k)]$ obtained by setting
$x_k=\exp (u_k+i\theta_k)$, the functions $g_f$ and $g_f'$ being defined on
 $\Rnm=\{-\infty\le u_k\le 0\}$.

\begin{prop}
 In order that a function $h(u_1,\ldots,u_n):\:\Rnm\to\Rm$  be the image of
$f\in PSH_-(D)$ obtained by $x_k=\exp (u_k+i\theta_k)$ and
 $$
f(x)=f[\exp (u_k+i\theta_k)]=h(u),
 $$
it is necessary and sufficient that $h$ be convex of $u\in\Rnm$, increasing
in each $u_k,\ -\infty\le u_k\le 0$, and $h(u)\not\equiv -\infty,\ h(u)\le
0$.
\label{prop:L2}
\end{prop}

The necessity condition results from the classic properties of
$f\in PSH_-(D)$. To show the sufficiency, we remark that convexity of $h$
implies its continuity on $\Rnm$. On the other hand, we have (the
derivatives being taken in the sense of distributions), for any
$\lambda\in\Cn$,
\beq
4\sum\frac{\partial^2}{\partial x_k \partial\bar x_j}\lambda_k
\overline{\lambda_j}
= \sum\frac{\partial^2}{\partial u_k \partial u_j}\lambda_k'
\overline{\lambda_j'}
\label{eq:L13}
\eeq
where $\lambda_k'=x_k^{-1}\lambda_k$. Let $A\subset D$ be the union of the
subspaces $\{x_k=0\}$ in $D$. By (\ref{eq:L13}),
$f\in PSH(D\setminus A)$. The condition $f(x)\le 0$ implies that
$f$  extends by upper semicontinuity to $A$, so $f\in PSH_-(D)$ for
 $f(x)=h(\log|x_1|,\ldots,\log|x_n|)$.

\medskip

{\bf Definition.} { We denote by $Conv\, (\Rnm)$ the class of
functions $h(u_1,\ldots,u_n)\le 0,\ -\infty\le u_k\le 0$, satisfying the
conditions listed in Proposition \ref{prop:L2}.

\medskip

\begin{prop}
Let $h\in Conv(\Rnm)$ be the image of
$f(x_k)=h[\exp (u_k+i\theta_k)]\in PSH_-(D)$. Then
\begin{enumerate}
\item[a)]
$$
\lim_{v\to -\infty}\,v^{-1}h(u_1+v,\ldots,u_n+v)=
 \lim_{v\to -\infty}\,{\partial\over \partial
v} h(u_1+v,\ldots,u_n+v)=\nu(f,0),
$$
where $\nu(f,0)$ is the Lelong number of $f$ at $x=0$.
More generally,
\beq
\lim_{v\to -\infty}\,v^{-1}h(u_1+a_1v,\ldots,u_n+a_nv)=
\nfoa
\label{eq:L14}
\eeq
is independent of $u_k$;
\item[b)]
$\lim_{v\to -\infty}\,v^{-1}h(u_1+v,u_2,\ldots,u_n)
  =\nu_1(f,0)$ is independent of $u_k$ and is the generic Lelong number
(cf. \cite{keyD4}) of the current $T= \dvapi dd^cf$ along the variety
$D_1=\{x\in D:\: x_1=0\}$. Moreover, for $x_1'=(x_2,\ldots,x_n)$,
 the function
$$
h_1(r_1,x_1')=(2\pi)^{-1}\int_0^{2\pi} f(r_1e^{i\theta_1},x_1')\,d\theta_1,
$$
has the property
$$
\lim_{w\to 0}\,(\log w)^{-1}h_1(wr_1,x_1')=\nu_1(f,0)
$$
for $x_1'\in D_1$  with exception of a pluripolar subset of $D_1$;
\item[c)]
 $\sum_1^n\,\nu_k(f,0)\le \nu(f,0)$.
\end{enumerate}
\label{prop:L3}
\end{prop}

{\it Proof.}
Existence and equality of the limits in $a)$  follow from the
increasing  with respect to $v,\ -\infty<v\le 0$, and from the condition
$h\le 0$. Moreover, if $l(\rho)$ is the mean value of $f(x)$ over the sphere
$|x|=\rho$, then
\beq
\nu(f,0)=\lim_{\rho\to 0}\,\frac{\partial\, l(\rho)}{\partial\log\rho}
 =\lim_{\rho\to 0}\,(\log\rho)^{-1}l(\rho).
\label{eq:L15}
\eeq
%
 We compare the mean values
with respect to $\theta_k$ over the circled
domains $B(0,\rho)$ and $D(\rho)=\{\sup |x_k|\le
\rho<1\}$ for the  image $h(u)$ of the circled
function $f_c(x)$,
 for
$u_k=\log \rho-{1\over 2}\log n$ and $u_k'=\log \rho,\ \okn$:
$$
h(u )\le l(\log\rho)\le h( u'),
$$
 since $D(\rho/\sqrt n)\subset B(0,\rho) \subset D(\rho)$.  This gives us
(\ref{eq:L14}) and $a)$.

Statement $b)$ is known (cf. \cite{keyL5}). The limit
\beq
-c(x_1')=\lim_{r\to 0}\, \left(\log\frac{1}{r}\right)^{-1}h(r,x_1')\le 0
\label{eq:L16}
\eeq
for $r\searrow 0$ exists and is obtained by increasing negative values,
the second term of (\ref{eq:L16}) belonging to $PSH(D_1)$ for $r>0$.
If $c(\hat x_1')=0$ for a point $\hat x_1' \in\cldo$, then  $c(x_1')=0$
except for a pluripolar subset of $D$ and the statement is proved.
Otherwise, consider the set $\overline{D_1}(r)\Subset\cldo$ and
$c_0=\sup c(x_1')$ for $x_1'\in D_1(r)$ and apply the preceding  argument
to
$$
\lim_{r\to 0}\, \left[\left(\log\frac{1}{r}\right)^{-1}h(r,x_1')+c_0\right].
$$
The statement for $h\in Conv(\Rnm)$ follows from this precise property
of the plurisubharmonic image.

To establish $c)$, we observe that for $u\in\Rnm$ and $h(u)$ the image in
$Conv(\Rnm)$ of $f\in PSH_-(\cld)$,
$$
{\partial\over \partial v} h(u_1+v,\ldots,u_n+v)=
\sum_1^n\,{\partial h\over \partial
u_k} (u_1+v,\ldots,u_n+v),
$$
the derivatives  are positive and decreasing for $v\searrow -\infty$, and
the limit of
$$
{\partial h\over \partial
u_k} (u_1+v,u_2,\ldots,u_n)
$$
 is equal to $\nu_1(f,0)$, the Lelong number of $dd^cf$ along $D_1$.
Therefore
$$
{\partial\over \partial v} h(u_1+v,\ldots,u_n+v)\ge \sum_1^n\nu_k(f,0),
$$
so taking $v\searrow -\infty$ we get by $a)$,
\beq
\nu(f,0)\ge \sum_1^n\nu_k(f,0).
\label{eq:L17}
\eeq

\medskip

{\it Remark.} Actually, by the theorem of Y.T. Siu, (\ref{eq:L17}) is a
particular case of the following statement: the number $\nu(f,0)$ is at
least equal to the sum of the generic numbers $\nu(W_i)$ for
$T= \dvapi dd^cf$ along analytic varieties $W_j$
of codimension $1$ containing the origin.

\medskip
In what follows, we will use a special subclass of circled
plurisubharmonic functions $f\in PSH_-(D)$ that have the following "conic"
property:  the convex image $g_f(u)$ of $f$ satisfies the equation
\beq
g_f(c\,u)=c\, g_f(u)\ {\rm for\ every}\ c>0.
\label{eq:conic}
\eeq
 Such a function $f$ will be called an {\it indicator}.  For example,
the weights $\vaa$ in  (\ref{eq:L5}) are indicators.

\begin{prop}
Let $f\in PSH_-(D)$ be an indicator.
Then
$(dd^cf)^n=0$ on $D_0=\{x\in D:\:x_1\ldots x_n\neq 0\}$.
\label{prop:conic}
\end{prop}

{\it Proof.} It is sufficient to show that  the domain $  D_0$
 can be foliated by one-dimensional analytic  varieties $\gamma_y $ such that
the restriction
 of $f$ to each leaf $\gamma_y$ is harmonic on $\gamma_y$.
So, given $y=(|y_k|e^{i\theta_k})\in D_0 $,  consider an analytic variety
$\gamma_y$, the image of $\bf C$ under the holomorphic mapping
$\lambda=(\lambda_1,\ldots,\lambda_n) $ with
$\lambda_k(\zeta)=|y_k|^\zeta e^{i\theta_k}$. Note that
$y=\lambda(1)\in\gamma_y$.
As  $f$ is circled, the function
$ f_y(\zeta)=f(\lambda(\zeta))$, the restriction
 of $f$ to   $\gamma_y$, depends only on ${\rm Re}\:\zeta$.
By (\ref{eq:conic}), $f_y(\zeta)$  satisfies
$f_y(c\,\zeta)=c\,  f_y(\zeta)$ for all $c>0 $.
Therefore, it is linear and thus harmonic on $\gamma_y$.

\section{Indicator of a plurisubharmonic function}

 Given a function $f\in PSH(\obl)$ and a point $x^0\in\obl$, we will
construct
 a function $\Ifxo(y)$ related to local properties of $f$ at $x^0$.
We will have $\Ifxo\in PSH_-(D)$, $D$ being the unit polydisk in the space
$\Cn_{(y)}$, and $ \Ifxo(y)<0$ in $D$ if and only if the Lelong number of $f$
at $x^0$ is strictly positive, otherwise $\Ifxo(y)\equiv 0$.

{\bf Definition.}  The {\it local indicator} (or just {\it indicator})
$\Ifo$  {\it  of a function $f\in
PSH_-(D)$, $D\subset\Cn_{(x)}$, at} $x^0=0$ is defined for $y\in
D\subset\Cn_{(y)}$ by
$$
\Ifo(y)=-n(f,0,-\log|y_k|).
$$
Referring to (\ref{eq:L9}) with $R=-\log w,\ 0<R<+\infty$, we rewrite this
as
\beq
\Ifo(y)=\lim_{R\to +\infty}\,R^{-1}f[\exp (u_k+i\theta_k+R\log|y_k|)].
\label{eq:L18}
\eeq
The limit (\ref{eq:L18}) exists {\it almost everywhere} for
$x_k=\exp (u_k+i\theta_k)$, however (see Introduction) the value
$\nfoa$  can be calculated as well by replacing $f(x)$ with
the circled functions $f_c(x)$ or $f_c'(x)$. One can then substitute them
for $f$ in (\ref{eq:L18}) to get $\Ifo$.
At $x^0\neq 0$, the function $\Ifxo(y)$ is defined by means of
$f[x_k^0+\exp (u_k+i\theta_k+R\log|y_k|)]$.

If $f$ is replaced by $f_c[\exp (u_k+i\theta_k)]=g_f(u_k)$ or by
$f_c'[\exp (u_k+i\theta_k)]=g_f'(u_k)$, the limit exists,
by Propositin \ref{prop:L3}, for {\it every}
$u=(u_k)\in\Rnm$:
\beq
\Ifo(y)=\lim_{R\to +\infty}\,R^{-1}g(u_k +R\log|y_k|)],
\label{eq:L19}
\eeq
  and does not depend on $u$.

\begin{prop} Let $f\in PSH_-(D)$. Then
\begin{enumerate}
\item[a)]
$\Ifo(y)\in PSH_-(D)$ and is $0$-circled;
\item[b)]
the convex image $g_\psi(u)$ in $\Rnm$ has the conic property
$ g_\psi(c\,u)=c\, g_\psi(u)$ for every $c>0$, i.e. $\Ifo$ is an indicator;
\item[c)]
\beq
\Ifo(y)\ge f_c'(y)\ge f(y),\ \forall y\in D;
\label{eq:L20}
\eeq
\item[d)]
the mapping $f\mapsto\Ifo$ is a projection, $\Ifo(y)$ is its  own
indicator at the origin;
\item[e)]
the indicator $\Ifo$ is the least indicator majorizing $f$ on $D$;
\item[f)]
if $f_j(x_j) $ is the restriction of $f$ to the complex
subspace $\{x_s=0,\ \forall s\neq j\}$
and
\beq
f_j(x_j)\not\equiv -\infty,
\label{eq:restr}
\eeq
 then
$\Ifo(y)\ge  \nu_j\log|y_j|$, $\nu_j$ being the Lelong number of
$f_j(x_j)$ at the origin;
\item[g)]
if (\ref{eq:restr}) holds for each $j$, then the Monge-Amp\`ere operator
$(dd^c\Ifo)^n$ is well defined on the whole polydisk $D$ and
\beq
(dd^c\Ifo)^n=0
\label{eq:MA}
\eeq
 on $D\setminus \{0\}$.
 \end{enumerate}
\label{prop:L4}
\end{prop}

{\it Proof.}
Statement $a)$ follows from (\ref{eq:L19}), $g(R\log|y_k|)$ being a
convex negative function for $R>0$, and the limit of the quotient is obtained
by increasing negative values. When setting $v_k=\log|y_k|=-a_k$, the image
of $\Ifo$ belongs to $Conv(\Rnm)$ and $\Ifo(y)$ is a $0$-circled
plurisubharmonic function.

The property $b)$, essential for the indicator
$\Ifo$, results from the  equality $n(f,0,c\, a)=c\,\nfoa$ for all $c>0$.

Relations $c)$ are a  consequence of (\ref{eq:L19}) where $g(u)$ is the
convex image $g_f'(u)$ of $f_c'(x)=\sup_{\theta_k}\,f(x_ke^{i\theta_k})$.
 We have $g_f'(\log|y_k|)\ge f(y)$. On the other hand, the quotient
$m(R)=R^{-1}g_f'(R\log|y_k|),\ R>R_0>1$,  is a convex, negative
and increasing function of $R$ for $|y_k|<1$. Therefore, $\lim_{R\to
+\infty}\,m(R)\ge m(1)$, and by (\ref{eq:L19}),
$$
0\ge\Ifo(y)\ge g_f'(\log|y_k|)\ge f(y)
$$
for $y\in D$.

Statement $d)$  follows from (\ref{eq:L19}) for $f=\Ifo$ and from
relation $b)$.

To prove $e)$, consider any indicator $\psi(y)\ge f(y)$ on $D$. Then
$\Psi_{\psi,0}(y)\ge \Ifo(y)$, and by $d)$, $\Psi_{\psi,0}=\psi$, so
$\psi(y)\ge\Ifo(y)\ \forall y\in D$.

The bound in $f)$ results from $c)$ and the maximum principle for
plurisubharmonic functions, since (for $j=1$)
$$
\Ifo(y)\ge \sup_{\theta_k}\, f(y_ke^{i\theta_k})\ge
\sup_{\theta_1}\, f(y_1e^{i\theta_1},0,\ldots,0)
$$
and for $|y_1|\searrow 0$ the quotient $(\log|y_1|)^{-1}\sup_{\theta_1}\,
f_1(y_1e^{i\theta_1})$ for the restriction $f_1$ to the complex
subspace $\{x_s=0,\ \forall s >1\}$, decreases to $\nu_1$.

Finally, in the  assumptions of $g)$, the function $\Ifo(y)$ is locally
bounded on $D\setminus \{0\}$ by $f)$, so the operator $(dd^c\Ifo)^n$ is well
defined on
$D$. Equation (\ref{eq:MA}) is valid on the domain
$D\setminus \{y:y_1y_2\ldots y_n=0\}$
 by Proposition \ref{prop:conic} and then on
$D\setminus \{0\}$, because the Monge-Amp\`ere measure of a bounded
plurisubharmonic function has zero mass on any pluripolar set
(see \cite{keyBT2}).

\medskip
{\it Remark}.  Statement $d)$ of Proposition \ref{prop:L4} is, in other
words, that  all the directional numbers $\nu(dd^c\Ifo, \vaa)$ of $\Ifo$
coincide  with the
directional numbers $\nu(dd^c f,\vaa)$ of the original function $f$, $\forall
a\in {\bf R}_+^n$.

\medskip

 The above construction is in fact of local character and
Proposition \ref{prop:L4} remains valid for the indicator
$\Ifxo$ of any function $f(x)$
plurisubharmonic in a neighbourhood $\omega$ of a  point $x^0\in\Cn$, with
the  the following change in the statement $c)$:   (\ref{eq:L20}) should be
replaced by
\begin{equation}
\Ifxo(x-x^0)\ge f(x) +  C\quad
\forall x\in D(x^0,r),\ C=C(u,r),
\label{eq:L20x}
\end{equation}
  where
$D(x^0,r)=\{ x:\:|x_k-x_k^0|< r,\ \okn\}$
and $r>0$ is such that the polydisk $D(x^0,r)\Subset\omega$.
And of course the restriction $f_j$ in (\ref{eq:restr}) should be taken
to the subspaces $\{x_s=x_s^0,\ \forall s\neq j\}$.

\medskip
 Let now $f(x)\in PSH(\omega)$ be locally bounded on $\omega\setminus
\{x^0\}$. Then its indicator $\Ifxo$ satisfies the equation
\beq
(dd^c\Ifxo)^n=\tfxo\,\delta(0)
\label{eq:tau1}
\eeq
with some number $\tfxo\ge 0$ and $\delta(0)$ the Dirac measure
at $0$, and $\tfxo> 0$ if and only if the Lelong number of the function
$f$ at $x^0$ is strictly positive.
And now we relate this value to $(dd^cf)^n$.

\medskip
\beth
Let $f\in PSH(\omega)$ be locally bounded out of a point $0\in\omega$. Then
\beq
(dd^cf)^n\ge\tfo\,\delta({0}).
\label{eq:ocenka}
\eeq
\label{theo:ocenka}
\eth

{\it Proof.}  In view of (\ref{eq:L20x}), the function $f$ satisfies
 \beq
\limsup_{x\to 0}\,\frac{\Ifo(x )}{f(x)}\le 1.
\label{eq:13}
\eeq
By the Comparison theorem of Demailly \cite{keyD4}, Theorem 5.9,
this implies
$$
(dd^c\Ifo)^n|_{\{0\}}\le (dd^cf)^n|_{\nol}.
$$
On the other hand,
$$
(dd^c\Ifo)^n|_{\nol}=(dd^c\Ifo)^n=\tfo\delta({0})
$$
by  (\ref{eq:tau1}),
that gives us (\ref{eq:ocenka}).

The theorem is proved.

\medskip
{\it Remark.}  It is well known  that for any plurisubharmonic function
$v$ with isolated singularity at $0$,  there is the relation
\beq
(dd^cv)^n\ge [\nu(dd^cv,0)]^n\delta({0}).
\label{eq:"1"}
\eeq
 By the remark  after the proof of Proposition \ref{prop:L4},
$\nu(dd^cf,0)$ is equal to $\nu(dd^c\Ifo,0)$. Applying (\ref{eq:"1"}) to
$v=\Ifo$ we get,
in view of
Theorem \ref{theo:ocenka},
$$
(dd^cf)^n\ge (dd^c\Ifo)^n\ge\left[\nu(dd^c\Ifo,0)\right]^n=
\left[\nu(dd^cf,0)\right]^n,
$$
so (\ref{eq:ocenka}) is an improvement of (\ref{eq:"1"}).

For example, if $f(x)= \log(|x_1|^{k_1}+|x_2|^{k_2})$ with
$0<k_1<k_2$, then
$$
(dd^cf)^2=\tau_f(0)\,\delta(0)=k_1k_2\,\delta(0)> k_1^2\,
\delta(0)=[\nu(dd^cf,0)]^2\delta(0),
$$
 and thus $\tau_f(0)>[\nu(dd^cf,0)]^2$.

More generally, if $F$ is a holomorphic mapping to $\Cn$ with an isolated
zero at $0$ of multiplicity $\mu_0$, and $f=\log|F|$, then
$$
[\nu ( dd^cf,0)]^n\le \tau_{ f}(0)\le \mu_0.
$$

\medskip
In fact, relation (\ref{eq:13}) makes it possible to obtain extra bounds for
$(dd^cf)^n$ in case of $exp\: f\in C(\obl)$. Such a function $f$ can be then
considered as a plurisubharmonic
weight $\varphi$ for Demailly's generalized numbers $\nu (T,\varphi)$ of a
closed positive current $T$ of bidimension $(p,p),\ 1\le p\le n-1$
\cite{keyD4}:
$$
\nu (T,\varphi)=\lim_{s\to -\infty}\int_{\{\varphi<s\}}T\wedge (dd^c
\varphi)^p=
T\wedge (dd^c\varphi)^p|_{\nol}.
$$
Moreover, the function $\Ifo$ is such a weight, too. By Comparison theorem
from \cite{keyD4}, Theorem 5.1, relation (\ref{eq:13}) implies
\beq
\nu(T,\Ifo )\le \nu(T,f).
\label{eq:14}
\eeq
Take
$$
T_k=(dd^cf)^k\wedge(dd^c\Ifo)^{n-k-1},\ \okn-1.
$$
These currents are well defined on a neighbourhood of $0$ and are of
bidimension $(1,1)$. Applying (\ref{eq:14}) to $T=T_k$ we obtain
$$
T_k\wedge dd^cf|_{\nol}\ge T_k\wedge dd^c\Ifo|_{\nol},
$$
that gives us

\medskip
\begin{prop}
Let $f\in PSH_-(\obl )$ be locally bounded out of $\nol$ and $\exp f\in
C(\obl)$. Then
\begin{eqnarray*}
(dd^cf)^n|_{\nol} &\ge &(dd^cf)^{n-1}\wedge dd^c\Ifo|_{\nol}\ge \ldots\\
&\ge &(dd^cf)^{n-k}\wedge (dd^c\Ifo)^k|_{\nol}\ge \ldots \ge (dd^c\Ifo)^n.
\end{eqnarray*}
\end{prop}

\section{Dirichlet problem with  local indicators }

Let $\obl$ be a bounded pseudoconvex domain in $\Cn$ and $K$ be a compact
subset of $\obl$. By $PSH(\obl,K)$ we denote the class of plurisubharmonic
functions on $\obl$ that are locally bounded on $\obl\setminus K$.

Let
$ K=\{x^1,\ldots,x^N\}\subset\obl$ and $\{\Psi_m\} $
be
$N$ indicators, i.e.  circled functions in $PSH_-(D)$ whose
convex images     satisfy (\ref{eq:conic}). In the sequel we  assume that
$\Ij \in PSH(D,\nol)$. Then by Proposition
\ref{prop:conic},
\beq
(dd^c\Ij)^n=\tau_m\delta(0),\ {\ojn}.
\label{eq:tau}
\eeq

 Let us fix the system $\Phi=\{(x^1,\Psi_1),\ldots,(x^N,\Psi_N)\}$ and
consider a positive measure $T_\Phi$ on $\obl$, defined
as
\beq
T_\Phi=\sum_{\ojn} \tau_m\delta({x^m}).
\label{eq:TPsi}
\eeq

Each function $\Ij$ can be extended from a neighbourhood of the origin to a
function $\tilde\Psi_m\in PSH(\Cn,\nol)$, and the indicators of the functions
\beq
\tilde\Psi(x)=\sum_m\tilde\Psi_m(x-x^m)+A,
\label{eq:summa}
\eeq
at $x^m$ are equal to $\Ij, \ojn$, for any real number $A$. So the class
\beq
N_{\Phi,\obl}=\{v\in PSH_-(\obl,K):\: \Psi_{v,x^m}\le\Ij,\,\ojn\}
\label{eq:class}
\eeq
is not empty.

Theorem \ref{theo:ocenka} implies

\medskip
\beth
$(dd^cf)^n\ge T_\Phi \quad\forall f\in N_{\Phi,\obl}.$
\label{theo:globoc}
\eth

\medskip
Now we introduce the function
\beq
 \green(x)=\sup\{v(x):\:v\in N_{\Phi,\obl}\}.
\label{eq:green}
\eeq

\medskip
\beth
Let $\obl$ be a hyperconvex domain in $\Cn$. Then the function
$ G=\green$ has the following properties:
\begin{enumerate}
\item[a)]
$ \gree\in PSH_-(\obl,K)$;
\item[b)]
$ \gree(x)\to 0$ as $x\to\partial\obl$;
\item[c)]
$\Psi_{ \gree,x^m}=\Ij,\ \ojn$;
\item[d)]
$(dd^c \gree)^n=T_{\Phi}$,   the measure
$T_{\Phi}$ being defined by (\ref{eq:TPsi});
\item[e)]
$ \gree\in C(\overline{\obl}\setminus K)$.
\end{enumerate}
\label{theo:green}
\eth

\medskip
{\it Remark}. In the case where $\Psi_m=\nu_m\log|x|$, the function
$ \green$, the pluricomplex Green function with several weighted
poles, was introduced in \cite{keyL3}. A situation with infinite number of
poles was considered in \cite{keyZ}, where a function
$G_{f,\obl}$ was introduced as the upper envelope of the class
$\{v\in PSH_-(\obl,K):\:\nu(dd^cv,x)\ge\nu(dd^cf,x),\ \forall x
\},\  f$ being a plurisubharmonic
function  with the following properties: $e^f$ is continuous,
$f^{-1}(-\infty)$ is a compact subset of $\obl$, and the set
$\{x:\nu(dd^cf,x)>0\}$ is dense in
$f^{-1}(-\infty)$. Our proof is much the same as of the corresponding
statements of \cite{keyZ}.

\medskip
{\it Proof of Theorem \ref{theo:green}.}
Since $N_{\Phi,\obl}\neq\emptyset$, the
function
$ G=\green$ is well defined and
$$
G^* =\limsup_{y\to x}\, \gree(y)\in PSH_-(\obl,K).
$$

The function $\tilde\Psi$ in (\ref{eq:summa}) can be modified
in a standard way to $\tilde\Psi'\in PSH_-(\obl,K)$ such that
$\tilde\Psi'(x)=\alpha\,\rho(x)$ in a neighbourhood of $\partial\obl$,
$\alpha$ being a positive number and $\rho$ a bounded exhaustion function on
$\obl$, and $\tilde\Psi'(x)=\tilde\Psi(x)-\beta$ on a neighbourhood of $K$.
It shows us that
\beq
G^* \ge \tilde\Psi'.
\label{eq:18}
\eeq
It implies, in particular, that
\beq
\Psi_{G^*,x^m}\ge \Psi_{\tilde\Psi',x^m}
=  \Ij,\ \ojn.
\label{eq:19}
\eeq

Since $\hpsi_{\sup\{v,w\},x} \le\sup\{\hpsi_{v,x},\hpsi_{w,x}\}$
for any plurisubharmonic functions $v$ and $w$,
 there exists an increasing sequence of functions $ v_j\in
N_{\Phi,\obl}$ such that $\lim_{j\to\infty}\, v_j= v\le\gree$
and $ v^*=\gree^*$.

 The indicator of $v_j$ at $x^m$
is   the limit of $R^{-1}g_{v_j,x^m}(R\log|y_k|)$ for $R\to
+\infty$, the function
$g_{v_j,x^m}(u)$ being the convex image of the mean value of
$v_j(x_k^m+e^{u_k+i\theta_k})$ with respect to $\theta_k $ for
$u_k<\log\: {\rm dist}\,(x^m,\partial\obl),\ \okn$, and
the limit is obtained by the increasing values. It gives us
\beq
R^{-1}g_{v_j,x^m}(R\log|y_k|)\le \Psi_{v_j,x^m}(y)\le\Ij(y).
\label{eq:20}
\eeq

 The functions $ v_j$ increase to $\gree^*$ out of a pluripolar set
$X=\{x\in\obl:\:  v(x)< v^*(x)\}$. Since the restriction of $X$ to
the distinguished boundary of any polydisk is of zero Lebesgue measure
\cite{keyL2}, (\ref{eq:20}) implies that
$$
R^{-1}g_{ \gree^*,x^m}(R\log|y_k|) \le\Ij(y).
$$
and thus, taking $R\to +\infty$,
\beq
\hpsi_{\gree^*,x^m}\le\Ij,\quad\ojn.
\label{eq:21}
\eeq
As $\gree^*\in PSH_-(\obl)$,  the function $\gree^*$ belongs to the class
$N_{\Phi,\obl}$ and so
\beq
\gree^*\equiv\gree.
\label{eq:rr}
\eeq

By (\ref{eq:19}) and (\ref{eq:21}), $\hpsi_{\gree,x^m}=\Ij$. It proves
statements $a)$ and $c)$; statement $b)$ follows from inequality
(\ref{eq:18}).

Continuity of $\gree$ can be  proved as in \cite{keyZ}, Theorem 2.6,
with the following modification. Instead of Demailly's approximation
theorem \cite{keyD5} we use the similar fact: for any function $u\in
PSH(\Omega)$ there exists a sequence of continuous plurisubharmonic
functions $u_m$ satisfying
$$
u(x)-{c_1\over m}\le u_m(x)\le\sup\,\{u(x+y):\: |y_k-x_k|\le r_k,
\ 1\le k\le n\}+{1\over m}\log{c_2\over r_1\ldots r_n}
$$
and
$$
\Psi_{u,x}(y)\le \Psi_{u_m,x}(y)\le \Psi_{u,x}(y)-{1\over m}
\log|y_1\ldots y_n|,\quad \forall x\in\Omega,\ \forall y\in D.
$$

To prove $d)$, observe that in view of (\ref{eq:18}) and (\ref{eq:rr}),
$\tilde\Psi'\le\gree$. By the comparison theorem of Demailly
(\cite{keyD4}, Theorem 5.9), this implies
$$
(dd^c\gree)^n|_{\{x^m\}}\le (dd^c \tilde\Psi')^n|_{\{x^m\}}\le
[dd^c\Ij(x-x^m)]^n,\ \ojn,
$$
and therefore
\beq
(dd^c\gree)^n|_{ K}\le\mera.
\label{eq:22}
\eeq
On the other hand, by Theorem \ref{theo:globoc},
$$
(dd^c\gree)^n\ge\mera.
$$
Being comparing to (\ref{eq:22}) this provides
$$
(dd^c\gree)^n|_{ K}=\mera.
$$
Finally, the equality $(dd^c\gree)^n=0$ on $\obl\setminus K$ can be proved
in a standard way by showing it is maximal on $\obl\setminus K$ (see
\cite{keyBT1}, \cite{keyD3}), that proves $d)$.

The theorem is proved.

\medskip
As a consequence, we get an "indicator" variant of the Schwarz type lemma
(see \cite{keyL3}, \cite{keyZ}):

\beth
Let the indicator   of a function $g\in PSH(\obl)$
at  $x^m $ does not exceed $\Ij$,   $ \ojn$,
  and let
$g(x)\le M$ on $\obl$. Then $g(x)\le M+\green (x),\ \forall x\in\obl$.
\label{theo:schwarz}
\eth

\medskip
Now we are going to show that the function $\green$ is the unique
plurisubharmonic
function with the properties $a)-d)$ of Theorem \ref{theo:green}.
It is known that for unbounded plurisubharmonic functions $u$,
the Dirichlet problem
\beq
 \cases{ (dd^cv)^n=\mu\ge 0 &
 on $ \obl$ \cr
 v=h &  on $\partial\obl $\cr}
\label{eq:diri}
\eeq
need not have a unique solution even in a simple case $\mu=\delta(0),\
h\equiv 0$. However, a solution is unique under some regularity
assumptions on the functions $v$. For example,
as was established in \cite{keyZ},   (\ref{eq:diri}) has a unique
solution
for
\beq
\mu=\sum[\nu(dd^cf,x^m)]^n\delta({x^m})
\label{eq:23}
\eeq
with $f(x)$  specified in the remark after the statement of Theorem
\ref{theo:green}, if the functions $v(x)\in PSH_-(\obl,K)$ have to satisfy
\beq
\nu(dd^cv,x^m)=\nu(dd^cf,x^m),\ \ojn.
\label{eq:24}
\eeq
These additional
relations mean that
\beq
v(x)\sim \nu(dd^cv,x^m)\,\log|x-x^m]\ {\rm near}\ x^m
\label{eq:reg}
\eeq
($v$ has regular densities at its poles, in the terminology of \cite{keyL3}).

In our situation,
\beq
\mu=T_f=\sum_m\tau_m\delta({x^m}),
\label{eq:Tf}
\eeq
where $\tau_m$ are defined by (\ref{eq:tau}) with $\Ij=\hpsi_{f,x^m}$,
and we are going to replace condition (\ref{eq:24})   by
$\Ivxm=\Ifxm,\
\ojn$.

To prove the uniqueness, we need a variant of the comparison theorem for
unbounded plurisubharmonic functions
(see \cite{keyBT1} - \cite{keyD1},
\cite{keyD3}, \cite{keyKl}, \cite{keyZ} for different classes of
plurisubharmonic functions).

\medskip
\beth
Let $f\in PSH(\obl,K),\ K=\{x^1,\ldots,x^m\}$, and
\beq
(dd^cf)^n|_{K}=T_f,
\label{eq:25}
 \eeq
the measure $T_f$ being given by (\ref{eq:Tf}).
Let $v\in PSH(\obl,K)$ satisfy the conditions
\begin{enumerate}
\item[1)]
$\liminf_{x\to\partial\obl}\,(f(x)-v(x))\ge 0$;
\item[2)]
$(dd^cv)^n\ge (dd^cf)^n $ on $\obl\setminus K$;
\item[3)]
$\Iuxm\le\Ifxm,\ \ojn$.
\end{enumerate}
Then $v\le f$ on $\obl$.
\label{theo:uniq}
\eth

\medskip
The proof is just as of Theorem 3.3 of \cite{keyZ}, and we omit it here.

\medskip

\begin{cor}
Under the conditions of Theorem \ref{theo:green}, the function $\green$
  is the unique plurisubharmonic function
with the properties $a)-d)$ of that theorem.
\end{cor}

\medskip

{\it Remarks.}

1.  Condition (\ref{eq:25}) is essential. Indeed, let
$$
f(x)={1\over 2}\log(|x_1|^4+|x_1+x_2^2|^2),\quad
v(x)={1\over 2}\log(|x_1|^2+|x_2|^4)+m
$$
 with
$m=\inf\,\{f(x):\:|x|=1\}>-\infty$. Then $v(x)\le f(x)$ for $|x|=1$,
$(dd^cf)^2=(dd^cv)^2=0$ on $\{0<|x|<1\}$,
and ${\Psi_{v,0}(x)}=\Ifo (x)=\log\,\max\,\{|x_1|, |x_2|^2\}$.
However, for $x_2=t\in (0,e^m),\ x_1=-x_2^2$, $f(x)=2\log t<\log t+m
<v(x)$. The reason here is that $(dd^cf)^2=4\delta(0)>2\delta(0)=T_f$.

2. By Comparison theorem of Demailly \cite{keyD4},  relation (\ref{eq:25})
is true when
$$
f(x) \sim \Ifxm (x-x^m) \ {\rm near}\ x^m,
$$
 a weaker than
(\ref{eq:reg}) but still controlled regularity.

3. By Proposition \ref{prop:L4}, an indicator $\hpsi$ possesses the
properties
$a)-d)$ of $\green$ from Theorem \ref{theo:green} with $\obl= D$, the
unit polydisk, and
$\Phi=(\{0\},\hpsi)$. Therefore, $\hpsi=G_{ D,\Phi}$.

\medskip

Theorems \ref{theo:green} and \ref{theo:uniq} allow us also to state
the following result.

\medskip
\beth
 Let $\obl$ be a bounded strictly pseudoconvex domain,
$K=\{x^1,\ldots,x^m\}$, and
let a function $f\in PSH(\obl,K) $ satisfy
 $$
(dd^cf)^n=T_f.
 $$
Then the Dirichlet
problem
$$
 \cases{ (dd^cv)^n= T_f &
 on $ \obl$\cr
\Iuxm=\Ifxm & for $\ojn$\cr
 v=h &  on $\partial\obl $\cr}
$$
has a unique solution in the class $PSH(\obl,K)$ for each function
$h\in C(\partial\obl)$. This solution is continuous on $\overline{\obl}
\setminus K$.
\label{theo:dir}
\eth

\medskip
{\it Proof.} Let $\Phi=\{(x^m,\Ij)\}$ with
 $\Ij=\Ifxm $.
Consider the class
$$
N_{f,h}=\{v\in PSH(\obl,K):\, \Ivxm\le\Ifxm\ \forall m, \ \lim_{x\to y}\,
v(x)=h(y)\ \forall y\in\partial\obl\}.
$$
Let $u_0(x)$ be the unique solution of the corresponding homogeneous problem
$$
 \cases{ (dd^cu)^n= 0 &
 on $ \obl$\cr
 u=h &  on $\partial\obl. $\cr}
$$
Then $u_0+\green\in N_{f,h}$, so $N_{f,h}\neq\emptyset$.

The desired solution $v_0$ is given as
$$
v_0(x)=\sup\,\{v(x):\:v\in N_{f,h}\}.
$$
Just as in the proof of Theorem \ref{theo:green}, one can show that $v_0$
does solve the problem and is continuous on $\overline\obl\setminus K$. The
uniqueness follows from Theorem \ref{theo:uniq}.

\medskip

Theorem \ref{theo:dir} can be related to the following question
wich was one of the motivations of the present study. Let
$F:\overline\obl\to\Cn$ be a holomorphic mapping with isolated zeros
$\{x^m\}\subset\obl$ of multiplicities $\mu_m$. Then the function
$f(x)= \log|F(x)|$ solves the Dirichlet problem
$$
 \cases{ (dd^cv)^n= \sum_m\,\mu_m\delta({x^m}) &
 on $ \obl$\cr
 v=f &  on $\partial\obl. $\cr}
$$
  Under what extra conditions on  $v$, the function $f$
is the unique solution of the problem?
By Theorem \ref{theo:dir}, if $f$ has regular behaviour at $x^m$ with respect
to its indicators, i.e. if
$$
(dd^c\Ifxm)^n=\mu_m\delta({x^m}),\ \ojn,
$$
it gives the unique solution to the problem
$$
 \cases{ (dd^cv)^n= \sum\mu_m\delta({x^m}) &
 on $ \obl$\cr
 \Iuxm=\Ifxm & $\ojn$\cr
 v=f &  on $\partial\obl. $\cr}
$$

  \par\bigskip
 {\small{\bf Acknowledgements.}   The second author is
grateful to the Insitut de Math\'ematiques de Jussieu  for the kind
hospitality.   }

\vskip 0.5cm

P. Lelong
\par
9, Place de Rungis
\par 75013 Paris, France

\vskip 0.5cm

A. Rashkovskii
\par
Mathematical Division, Institute for Low Temperature Physics
\par
47 Lenin Ave., Kharkov 310164,
Ukraine

\vskip0.1cm

E-mail: \quad rashkovskii@ilt.kharkov.ua \quad rashkovs@ilt.kharkov.ua



\end{document}